\documentclass[oneside,english]{amsart}
\usepackage{amsthm}
\usepackage{amstext}
\usepackage{setspace}
\usepackage{amssymb}

\newcommand{\noun}[1]{\textsc{#1}}
\providecommand{\tabularnewline}{\\}

\numberwithin{equation}{section}
\numberwithin{figure}{section}
\theoremstyle{plain}
\newtheorem{thm}{Theorem}
  \theoremstyle{plain}
  \newtheorem{lem}[thm]{Lemma}
  \newcounter{casectr}
  \newenvironment{caseenv}
  {\begin{list}{{\itshape\ Case} \arabic{casectr}.}{%
   \setlength{\leftmargin}{\labelwidth}
   \addtolength{\leftmargin}{\parskip}
   \setlength{\itemindent}{\listparindent}
   \setlength{\itemsep}{\medskipamount}
   \setlength{\topsep}{\itemsep}}
   \setcounter{casectr}{0}
   \usecounter{casectr}}
  {\end{list}}
  \theoremstyle{plain}
  \newtheorem*{thm*}{Theorem}

\usepackage{amsmath}

\date{\today}

\usepackage{babel}

\begin{document}

\title[Coefficients of Level $p$ Modular Functions]{Divisibility Properties of
Coefficients of Level $p$ Modular Functions for Genus Zero Primes}

\author{Nickolas Andersen}
\address{Brigham Young University Department of Mathematics, Provo, UT 84602}
\curraddr{} \email{nickolasandersen@gmail.com}
\thanks{}

\author{Paul Jenkins}
\address{Brigham Young University Department of Mathematics, Provo, UT 84602}
\curraddr{} \email{jenkins@math.byu.edu}
\thanks{}

%\subjclass[2010]{Primary 11F33}
\begin{abstract}
Lehner's 1949 results on the $j$-invariant showed high divisibility
of the function's coefficients by the primes $p\in\{2,3,5,7\}$. Expanding
his results, we examine a canonical basis for the space of level $p$
modular functions holomorphic at the cusp $0$. We show that the Fourier
coefficients of these functions are often highly divisible by these
same primes.
\end{abstract}
\maketitle

\section{Introduction and Statement of Results}

A level $p$ modular function $f(\tau)$ is a holomorphic function
on the complex upper half-plane which satisfies \[
f\left(\frac{a\tau+b}{c\tau+d}\right)=f(\tau)\text{ for all }\left(\begin{smallmatrix}a & b\\
c & d\end{smallmatrix}\right)\in\Gamma_{0}(p)\]
 and is meromorphic at the cusps of $\Gamma_{0}(p)$. Equivalently,
$f(\tau)$ is a weakly holomorphic modular form of weight $0$ on
$\Gamma_{0}(p)$. Such a function will necessarily have a $q$-expansion
of the form $f(\tau)=\sum_{n=n_{0}}^{\infty}a(n)q^{n}$, where $q=e^{2\pi i\tau}$.

Of particular interest in the study of modular forms is the
classical $j$-invariant, $j(\tau)=q^{-1}+744 + \sum_{n=1}^{\infty}
c(n)q^{n}$, which is a modular function of level $1$. The
coefficients $c(n)$ of the $j$-function, like the Fourier
coefficients of many other modular forms, are of independent
arithmetic interest; for instance, they appear as dimensions of a
special graded representation of the Monster group.

In 1949 Lehner showed \cite{Lehner:1,Lehner:2} that \[
c(2^{a}3^{b}5^{c}7^{d}n)\equiv0\pmod{2^{3a+8}3^{2b+3}5^{c+1}7^{d}},\]
 proving that the coefficients $c(n)$ are often highly divisible
by small primes. Similar results have recently been proven for other
modular functions in~\cite{Griffin}, and for modular forms of level
1 and small weight in~\cite{Doud:padic}, \cite{DoudJenkinsLopez}. It
is natural to ask whether such congruences hold for the Fourier
coefficients of modular functions of higher level, such as those
studied by Ahlgren~\cite{Ahlgren:theta} in his work on Ramanujan's
$\theta$-operator.

Lehner's results for $j(\tau)$ are in fact more general; in~\cite{Lehner:2}
he pointed out that for $p=2,3,5,7$, similar congruences hold for
the coefficients of level $p$ modular functions which have integral
coefficients at both cusps and have poles of order less than $p$
at the cusp at infinity.

In this paper, for $p\in\{2,3,5,7\}$, we examine canonical bases for
spaces of level $p$ modular functions which are holomorphic at the
cusp $0$. To construct these bases, we introduce the level $p$
modular function $\psi^{(p)}(\tau)$, defined as \[
\psi^{(p)}(\tau)=\left(\frac{\eta(\tau)}{\eta(p\tau)}\right)^{\frac{24}{p-1}}
\text{ where }\eta(\tau) =
q^{\frac{1}{24}}\prod_{n=1}^{\infty}(1-q^{n}).\] The integer
$\frac{24}{p-1}$ for $p=2,3,5,7$ will appear frequently, so we will
denote it $\lambda^{(p)}$, or simply $\lambda$ where no confusion
arises. The function $\psi^{(p)}(\tau)$ is a modular function of
level $p$ with a simple pole at $\infty$ and a simple zero at 0.  We
will also use the modular function \[\phi^{(p)}(\tau) =
(\psi^{(p)}(\tau))^{-1}.\]

Following Ahlgren \cite{Ahlgren:theta}, and using the notation of
Duke and Jenkins \cite{Duke:zeros}, for $p=2,3,5,7$ we construct a
basis $\{f_{0,m}^{(p)}(\tau)\}_{m=0}^{\infty}$ for the space of
level $p$ modular functions which are holomorphic at 0 as follows:
\[ f_{0,0}^{(p)}(\tau)=1,\]
\[f_{0,m}^{(p)}(\tau)=q^{-m}+O(1) =
\psi^{(p)}(\tau)^{m}-Q(\psi^{(p)}(\tau)),\] where $Q(x)$ is a
polynomial of degree $m-1$ with no constant term, chosen to
eliminate all negative powers of $q$ in $\psi^{(p)}(\tau)^{m}$
except for $q^{-m}$. Since $\psi^{(p)}(\tau)$ vanishes at $0$ and
the polynomial $Q$ has no constant term, we see that the functions
$f_{0,m}^{(p)}$ also vanish at $0$ when $m>0$. We write\[
f_{0,m}^{(p)}=q^{-m}+\sum_{n=0}^{\infty}a_{0}^{(p)}(m,n)q^{n}\]
 so that for $n\geq0$, the symbol $a_{0}^{(p)}(m,n)$ denotes the coefficient
of $q^{n}$ in the $m^{th}$ basis element of level $p$.  Note that
the function $f_{0, m}^{(p)}$ corresponds to Ahlgren's $j_m^{(p)}$.

For an example of some of these functions, consider the case $p=2$:
\begin{align*}
f_{0,1}^{(2)}(\tau) & =\psi^{(2)}(\tau)\\
 & =q^{-1}-24+276q-2048q^{2}+11202q^{3}-49152q^{4}+\ldots\\
f_{0,2}^{(2)}(\tau) & =\psi^{(2)}(\tau)^{2}+48\psi^{(2)}(\tau)\\
 & =q^{-2}-24-4096q+98580q^{2}-1228800q^{3}+10745856q^{4}+\ldots\\
f_{0,3}^{(2)}(\tau) & =\psi^{(2)}(\tau)^{3}+72\psi^{(2)}(\tau)^{2}+900\psi^{(2)}(\tau)\\
 & =q^{-3}-96+33606q-1843200q^{2}+43434816q^{3}-648216576q^{4}+\ldots\end{align*}
The function $f_{0,m}^{(p)}$ is a level $p$ modular function that
vanishes at $0$ (if $m\neq0$) and has a pole of order $m$ at
$\infty$. The conditions at the cusps determine this function
uniquely; if two such functions exist, their difference is a
holomorphic modular function, which must be a constant. Since both
functions vanish at 0, this constant must be $0$.

The functions comprising these bases for $p=2,3,5,7$ have
divisibility properties which bear a striking resemblance to the
divisibility properties of $j(\tau)$; in many cases they are
identical. As an example of some of the divisibility properties we
encounter with this basis, we experimentally examine the $2$-adic
valuation of the even indexed coefficients of $f_{0,m}^{(2)}(\tau)$
for $m=1,3,5,7$ in Table \ref{tab:2-Adic-Table}.  As the data in the
table suggest, the $2$-divisibility which $j(\tau)$ exhibits gives
us a lower bound on the $2$-divisibility of the odd-indexed $p=2$
basis elements.

\begin{table}[h]
\noindent \begin{centering}
\label{Flo:2-Adic-Float}
\par\end{centering}
\noindent \begin{centering}
\begin{tabular}{|c|c|c|c|c|c|c|c|}
\hline
\multicolumn{1}{|c}{} &  & \multicolumn{1}{c}{$a_{0}^{(2)}(m,2)$} & \multicolumn{1}{c}{$a_{0}^{(2)}(m,4)$} & \multicolumn{1}{c}{$a_{0}^{(2)}(m,6)$} & \multicolumn{1}{c}{$a_{0}^{(2)}(m,8)$} & \multicolumn{1}{c}{$a_{0}^{(2)}(m,10)$} & $a_{0}^{(2)}(m,12)$\tabularnewline
\hline
 & $1$ & $11$ & $14$ & $13$ & $17$ & $12$ & $16$\tabularnewline
\cline{2-8}
 & $3$ & $13$ & $16$ & $15$ & $19$ & $14$ & $18$\tabularnewline
\cline{2-8}
$m$ & $5$ & $12$ & $15$ & $14$ & $18$ & $13$ & $17$\tabularnewline
\cline{2-8}
 & $7$ & $14$ & $17$ & $16$ & $20$ & $15$ & $19$\tabularnewline
\cline{2-8}
 & min & $11$ & $14$ & $13$ & $17$ & $12$ & $16$\tabularnewline
\hline
\multicolumn{8}{|c|}{}\tabularnewline
\hline
$j(\tau)$ &  & $11$ & $14$ & $13$ & $17$ & $12$ & $16$\tabularnewline
\hline
\end{tabular}
\par\end{centering}
\noindent \centering{}\caption{\label{tab:2-Adic-Table}$2$-adic
valuation of $a_{0}^{(2)}(m,n)$ compared to corresponding
coefficients in $j(\tau)$}
\end{table}

Note that these functions form a basis for $M_0^\infty(p)$, the
space of modular forms of weight 0 and level $p$ with poles allowed
only at the cusp at $\infty$.  A full basis for the space $M_0^!(p)$
of weakly holomorphic modular forms of weight $0$ and level $p$ is
generated by the $f_{0, m}^{(p)}(\tau)$ and the functions
$(\phi^{(p)}(\tau))^n$ for integers $n\geq 1$.

Recall that the concluding remarks of Lehner's second paper
\cite{Lehner:2} state that the coefficients of certain level $p$
modular functions having a pole of order less than $p$ at $\infty$
have the same $p$-divisibility properties as the coefficients $c(n)$
of $j(\tau)$. More precisely, we have the following theorem.
\begin{thm}[Lehner]
\label{thm:Lehner-Main} Let $p\in\{2,3,5,7\}$ and let $f(\tau)$ be a
modular function on $\Gamma_{0}(p)$ having a pole at $\infty$ of
order $<p$ and $q$-series of the form
\[f(\tau)=\sum_{n=n_{0}}^{\infty}a(n)q^{n},\]
\[ f(-1/p\tau) = \sum_{m=m_0}^\infty b(n) q^n,\]
where $a(n), b(n) \in\mathbb{Z}$. Then the coefficients $a(n)$
satisfy the following congruence properties:
\[
\begin{array}{rll}
a(2^{a}n)\equiv0 & \pmod{2^{3a+8}} & \text{if }p=2\\
a(3^{a}n)\equiv0 & \pmod{3^{2a+3}} & \text{if }p=3\\
a(5^{a}n)\equiv0 & \pmod{5^{a+1}} & \text{if }p=5\\
a(7^{a}n)\equiv0 & \pmod{7^{a}} & \text{if }p=7.\end{array}\]
\end{thm}
Note that Lehner's original statement of this theorem mistakenly
states that a function on $\Gamma_{0}(p)$ inherits the
$p$-divisibility property for \emph{every }prime in $\{2,3,5,7\}$,
not just the prime matching the level.

A necessary condition in the statement of Lehner's theorem is that
the function must have an integral $q$-expansion at $0$. This
condition is quite strong; in fact, neither the function
$\phi^{(p)}(\tau)$ nor any of its powers satisfy it, although the
functions $f_{0, m}^{(p)}(\tau)$ do.

Further, Lehner's theorem assumes that the order of the pole at
$\infty$ must be less than $p$. In this paper, we remove this
restriction on the order of the pole to show that every function in
the $f_{0,m}^{(p)}$ basis has divisibility properties similar to
those in Theorem \ref{thm:Lehner-Main}.  Specifically, we prove the
following theorem.
\begin{thm}
\label{thm:Andersen-Main}Let $p\in\{2,3,5,7\}$, and let \[
f_{0,m}^{(p)}(\tau)=q^{-m}+\sum\limits _{n=0}^{\infty}a_{0}^{(p)}(m,n)q^{n}\]
 be an element of the basis described above, with $m=p^{\alpha}m'$
and $(m',p)=1$. Then, for $\beta > \alpha$,
\[
\begin{array}{rll}
a_{0}^{(2)}(2^{\alpha}m',2^{\beta}n)\equiv0 & \pmod{2^{3(\beta-\alpha)+8}} & \text{if }p=2\\
a_{0}^{(3)}(3^{\alpha}m',3^{\beta}n)\equiv0 & \pmod{3^{2(\beta-\alpha)+3}} & \text{if }p=3\\
a_{0}^{(5)}(5^{\alpha}m',5^{\beta}n)\equiv0 & \pmod{5^{(\beta-\alpha)+1}} & \text{if }p=5\\
a_{0}^{(7)}(7^{\alpha}m',7^{\beta}n)\equiv0 & \pmod{7^{(\beta-\alpha)}} & \text{if }p=7.\end{array}\]
\end{thm}

Note that for basis elements $f_{0,m}^{(p)}$ with $(m,p)=1$, these
divisibility properties match those in Theorem
\ref{thm:Lehner-Main}; in fact, Lehner's proof is easily extended to
prove the congruences in these cases. For basis elements with
$m=p^{\alpha}m'$ and $\alpha\ge1$, the divisibility is ``shifted.''
This shifting occurs in the $(\beta-\alpha)$ factor in the exponent
of the modulus.

For the coefficients $a_0^{(p)}(p^\alpha m', p^\beta n)$ with
$\alpha > \beta$, computations suggest that similar congruences
should hold.  Additionally, it appears that powers of the function
$\phi^{(p)}(\tau)$ have Fourier coefficients with slightly weaker
divisibility properties, despite the fact that their Fourier
coefficients at $0$ are not integral. It would be interesting to
more fully understand these congruences.

\section{Preliminary Lemmas and Definitions}

In this section, we provide the necessary definitions and background
for the proof of the main theorem.

For a prime $p$ we define the level $p$ Hecke operator $U_{p}$ by\[
U_{p}f(\tau)=\frac{1}{p}\sum_{\ell=0}^{p-1}f\left(\frac{\tau+\ell}{p}\right),\]
using the notation $U_{p}^{n}f=U_{p}U_{p}\cdots U_{p}f$ for repeated
application of $U_{p}$. When $f$ has the Fourier expansion
$f(\tau)=\sum_{n=n_{0}}^{\infty}a(n)q^{n}$, this operator takes the
form\[ U_{p}f(\tau)=\sum_{n=n_{0}}^{\infty}a(pn)q^{n},\] essentially
{}``pulling out'' all of the coefficients of $f$ whose index is
divisible by $p$. This operator preserves modularity:
%in fact,
if $f$ is a level $p$ modular function, then $U_{p}f$ is also a level
$p$ modular function.

For the primes $p=2,3,5,7$ the topological genus of
$\Gamma_{0}(p)\backslash\mathcal{H}$ is zero, so the field of level
$p$ modular functions is generated by a single modular function
called a Hauptmodul.  For the primes in consideration, one such
function is $\psi^{(p)}(\tau)$. Note that the modular function
$\phi^{(p)}(\tau)=\psi^{(p)}(\tau)^{-1}=q+O(q^{2})$ is also a
Hauptmodul.

Further, for these primes, the fundamental domain for
$\Gamma_{0}(p)$ has precisely two cusps, which may be taken to be at
$\infty$ and at $0$. Hence, we are most concerned with the behavior
of these functions at $\infty$ and at $0$. To switch between cusps,
we make the substitution $\tau\mapsto-1/p\tau$. The following lemma
gives relations for $\psi^{(p)}(\tau)$ and $\phi^{(p)}(\tau)$ at
$0$, and makes clear that powers of $\phi^{(p)}$ do not satisfy
Lehner's integrality condition.
\begin{lem}
\label{lem:Psi-At-0}The functions $\psi^{(p)}(\tau)$ and $\phi^{(p)}(\tau)$
satisfy the relations \begin{align}
\psi^{(p)}(-1/p\tau) & =p^{\lambda/2}\phi^{(p)}(\tau)\label{eq:Psi-At-0}\\
\phi^{(p)}(-1/p\tau) & =p^{-\lambda/2}\psi^{(p)}(\tau)\label{eq:Phi-At-0}\end{align}
\end{lem}
\begin{proof}
The functional equation for $\eta(\tau)$ is
$\eta(-1/\tau)=\sqrt{-i\tau}\eta(\tau)$. Using this, we find that \[
\psi^{(p)}\left(\frac{-1}{p\tau}\right) =
\left(\frac{\eta(-1/p\tau)} {\eta(-1/\tau)}\right)^{\lambda} =
\left(\frac{\sqrt{-ip\tau}
 \eta(p\tau)}{\sqrt{-i\tau} \eta(\tau)}\right)^{\lambda} = (\sqrt{p})^{\lambda}
\left(\frac{\eta(p\tau)}{\eta(\tau)}\right)^{\lambda} =
p^{\lambda/2}\phi^{(p)}(\tau).\] The second statement follows after
replacing $\tau$ by $-1/p\tau$ in the first statement.
\end{proof}
We next state a well-known lemma which gives a formula for
determining the behavior of a modular function at $0$ after $U_{p}$
has been applied. A proof can be found in \cite{Apostol:modular}.
\begin{lem}
\label{lem:Main-Up-Formula}Let $p$ be prime and let $f(\tau)$ be
a level $p$ modular function. Then\begin{equation}
p(U_{p}f)(-1/p\tau)=p(U_{p}f)(p\tau)+f(-1/p^{2}\tau)-f(\tau).\label{eq:Main-Up-Formula}\end{equation}

\end{lem}
Lehner's original papers included the following lemma and its proof, which gives
two important equations. The first gives a formula for $U_{p}\phi^{(p)}$
as a polynomial with integral coefficients in $\phi^{(p)}$; the second
gives an algebraic relation which is satisfied by $\phi^{(p)}(\tau/p)$.
\begin{lem}
\label{lem:Modular-Eq-Phi}Let $p\in\{2,3,5,7\}$. Then there exist
integers $b_{j}^{(p)}$ such that

\[
\begin{array}{cc}
\text{(a)} & U_{p}\phi^{(p)}(\tau)=p\sum\limits _{j=1}^{p}b_{j}^{(p)}\phi^{(p)}(\tau)^{j}.\end{array}\]

Further, let $h^{(p)}(\tau)=p^{\lambda/2}\phi^{(p)}(\tau/p).$ Then\[
\begin{array}{cc}
\text{(b)} & \big(h^{(p)}(\tau)\big)^{p}+\sum\limits _{j=1}^{p}(-1)^{j}g_{j}(\tau)\big(h^{(p)}(\tau)\big)^{p-j}=0\end{array}\]

where $g_{j}(\tau)=(-1)^{j+1}p^{\lambda/2+2}\sum\limits _{\ell=j}^{p}b_{\ell}\phi^{(p)}(\tau)^{\ell-j+1}$.\end{lem}
\begin{proof}

(a) Since $\phi$ vanishes at $\infty$, $U_{p}\phi$ also vanishes
at $\infty$; we will now consider its behavior at $0$. Using (\ref{eq:Main-Up-Formula})
and replacing $\tau$ by $p\tau$ in (\ref{eq:Psi-At-0}) we obtain\begin{align*}
U_{p}\phi(-1/p\tau) & =U_{p}\phi(p\tau)+p^{-1}\phi(-1/p^{2}\tau)-p^{-1}\phi(\tau)\\
 & =U_{p}\phi(p\tau)+p^{-1}\psi(p\tau)-p^{-1}\phi(\tau)\\
 & =O(q^{p})+p^{-\lambda/2-1}q^{-p}+O(1)-p^{-1}q+O(q^{2})\\
p^{\lambda/2+1}U_{p}\phi(-1/p\tau) & =q^{-p}+O(1)\end{align*}

The right side of this equation is a level $p$ modular function with
integer coefficients, so we may write it as a polynomial in
$\psi(\tau)$ with integer coefficients. The polynomial will not have
a constant term since the left side vanishes at $0$. Therefore,\[
p^{\lambda/2+1}U_{p}\phi(-1/p\tau)=\sum_{j=1}^{p}c_{j}\psi(\tau){}^{j}.\]

Now, replacing $\tau$ by $-1/p\tau$, we find\[
p^{\lambda/2+1}U_{p}\phi(\tau)=\sum_{j=1}^{p}c_{j}p^{\lambda j/2}\phi(\tau){}^{j}.\]

After cancelling the $p^{\lambda/2+1}$, we find that
$U_{p}\phi(\tau)=\sum\limits _{j=1}^{p}c_{j}'\phi(\tau)^{j}$ and we
compute the coefficients $c_{j}'$ (the authors used
\noun{mathematica}). The computation is finite, and we find that
each coefficient $c_{j}'$ has a factor of $p$, so the coefficients
$b_{j}^{(p)}$ are integral. A complete table of values of the
$b_{j}^{(p)}$ is found in Table \ref{tab:b_j-Table}.

\begin{table}[h]
\noindent \begin{centering}
\label{Flo:b_j-Float}
\par\end{centering}

\noindent \begin{centering}
\begin{tabular}{|c|c|c|c|c|c|}
\cline{3-6}
\multicolumn{1}{c}{} &  & \multicolumn{4}{c|}{$p$}\tabularnewline
\cline{3-6}
\multicolumn{1}{c}{} &  & 2 & 3 & 5 & 7\tabularnewline
\hline
 & 1 & $3\cdot2^{2}$ & $10\cdot3^{1}$ & $63\cdot5^{0}$ & $82\cdot7^{0}$\tabularnewline
\cline{2-6}
 & 2 & $2^{10}$ & $4\cdot3^{6}$ & $52\cdot5^{3}$ & $176\cdot7^{2}$\tabularnewline
\cline{2-6}
 & 3 &  & $3^{10}$ & $63\cdot5^{5}$ & $845\cdot7^{3}$\tabularnewline
\cline{2-6}
$j$ & 4 &  &  & $6\cdot5^{8}$ & $272\cdot7^{5}$\tabularnewline
\cline{2-6}
 & 5 &  &  & $5^{10}$ & $46\cdot7^{7}$\tabularnewline
\cline{2-6}
 & 6 &  &  &  & $4\cdot7^{9}$\tabularnewline
\cline{2-6}
 & 7 &  &  &  & $7^{10}$\tabularnewline
\hline
\end{tabular}
\par\end{centering}

\noindent \centering{}\caption{\label{tab:b_j-Table}Values of $b_{j}^{(p)}$}

\end{table}

(b) We again apply (\ref{eq:Main-Up-Formula}) to $\phi(\tau)$, this
time using what we know from $(a)$.\[
pU_{p}\phi(-1/p\tau)=pU_{p}\phi(p\tau)+\phi(-1/p^{2}\tau)-\phi(\tau)\]
\[
p^{2}\sum_{j=1}^{p}b_{j}\phi(-1/p\tau)^{j}=p^{2}\sum_{j=1}^{p}b_{j}\phi(p\tau)^{j}+\phi(-1/p^{2}\tau)-\phi(\tau).\]
We now use Lemma \ref{lem:Psi-At-0} with the knowledge that $\psi(\tau)=\phi(\tau)^{-1}$
to obtain
\[p^{2}\sum_{j=1}^{p}b_{j}p^{-\lambda j/2}\phi(\tau)^{-j}-p^{2}\sum_{j=1}^{p}b_{j}\phi(p\tau)^{j}+\phi(\tau)-p^{-\lambda/2}\phi(p\tau)^{-1}=0.\]
After replacing $\tau$ by $\tau/p$ and multiplying by $p^{\lambda/2}$,
this becomes\begin{equation}
p^{\lambda/2+2}\sum_{j=1}^{p}b_{j}\big(h(\tau)^{-j}-\phi(\tau)^{j}\big)+h(\tau)-\phi(\tau)^{-1}=0.\label{eq:Intermediate-Modular-Eq-Phi}\end{equation}

We now divide by $h^{-1}-\phi$. Note two facts: \[
h^{-j}-\phi^{j}=(h^{-1}-\phi)\sum_{\ell=0}^{j-1}h^{-\ell}\phi^{j-\ell-1}\]
\[
\frac{h-\phi^{-1}}{h^{-1}-\phi}=\frac{h(h\phi-1)}{\phi(1-h\phi)}=-\frac{h}{\phi}.\]
So (\ref{eq:Intermediate-Modular-Eq-Phi}) becomes\[
p^{\lambda/2+2}\sum_{j=1}^{p}b_{j}\sum_{\ell=0}^{j-1}h^{-\ell}\phi^{j-\ell-1}-\phi^{-1}h=0\]
which, after multiplying by $\phi h^{p-1}$, becomes\[
p^{\lambda/2+2}\sum_{j=1}^{p}b_{j}\sum_{\ell=0}^{j-1}h^{p-\ell-1}\phi^{j-\ell}-h^{p}=0.\]
We now change the order of summation and rearrange terms to obtain
the desired formula:\[
h(\tau)^{p}=\sum_{j=1}^{p}\big(p^{\lambda/2+2}\sum_{\ell=j}^{p}b_{\ell}\phi(\tau)^{\ell-j+1}\big)h(\tau){}^{p-j}.\]

\end{proof}
The next lemma states that when you apply $U_{p}$ to a certain type
of polynomial in $\phi_{p}$, you get a similar polynomial back which
has picked up a power of $p$. The details of this lemma are found in
both \cite{Lehner:1} and \cite{Lehner:2}, scattered throughout the
proofs of the main theorems. For our purposes, it will be more
useful in the following form.
\begin{lem}
\label{lem:Phi-Polynomials}Let $p\in\{2,3,5,7\}$ and let $R^{(p)}$
denote the set of polynomials in $\phi^{(p)}$ of the form\[
d_{1}\phi^{(p)}(\tau)+\sum_{n=2}^{N}d_{n}p^{\gamma}\phi^{(p)}(\tau)^{n}\]
\[
\begin{array}{ll}
\text{where }\gamma= & \begin{cases}
8(n-1) & \text{if }p=2\\
4(n-1) & \text{if }p=3\\
n & \text{if }p=5\\
n & \text{if }p=7.\end{cases}\end{array}\] Then $U_{p}$ maps
$R^{(p)}$ to $p^{\delta}R^{(p)}$ where $\delta=3,2,1,1$ for
$p=2,3,5,7$, respectively. That is, applying $U_{p}$ to a polynomial
of the above form yields a polynomial of the same form with an extra
factor of $p^{\delta}$.\end{lem}
\begin{proof}
Consider the function\[
d_{1}U_{p}\phi(\tau)+\sum_{n=2}^{r}d_{n}p^{\gamma}U_{p}\phi(\tau)^{n}.\]
For the first term, Lemma \ref{lem:Modular-Eq-Phi}(a) shows that
$U_{p}\phi(\tau)\in p^{\delta}R_{p}$ since, by inspection, the $b_{j}^{(p)}$
integers are divisible by sufficiently high powers of $p$.

For the remaining terms, we will prove\begin{equation}
p^{\gamma}U_{p}\phi^{n}=p^{\delta}r\label{eq:Up-Phi^t}\end{equation}
where $r\in R_{p}$. The result will immediately follow.

By the definition of $U_{p}$ we have\begin{equation}
U_{p}\phi^{n}=p^{-1}\sum_{\ell=0}^{p-1}\phi\left(\frac{\tau+\ell}{p}\right)^{n}=p^{-1-\lambda t/2}\sum_{\ell=0}^{p-1}h_{\ell}(\tau)^{n}\label{eq:Up-Def-W-Sum}\end{equation}
where $h_{\ell}(\tau)=p^{\lambda/2}\phi\left(\frac{\tau+\ell}{p}\right)$
is related to $h$ from Lemma \ref{lem:Modular-Eq-Phi}(b). Let $S_{n}$
be the sum of the $n^{th}$ powers of the $h_{\ell}$ so that \[
S_{n}=\sum_{\ell=0}^{p-1}h_{\ell}^{n}.\]

Define the polynomial
$F(x)=\sum_{j=0}^{p}(-1)^{j}g_{j}(\tau)x^{p-j}$ where the
$g_{j}(\tau)$ are as in Lemma \ref{lem:Modular-Eq-Phi}. In the same
lemma, if we replace $\tau$ with $\tau+\ell$, the $g_{j}(\tau)$ are
unaffected since $\phi(\tau+1)=\phi(\tau)$. Therefore, that lemma
tells us that the $p$ roots of the polynomial $F(x)$ are precisely
the $h_{\ell}$. Using Newton's formula for the $n^{th}$ power sum of
the roots of a polynomial, we obtain\begin{equation}
S_{n}=\sum_{\ell=0}^{p-1}h_{\ell}^{n}=\sum_{j=1}^{n}(-1)^{j+1}g_{j}S_{n-j}\label{eq:Newtons-Formula}\end{equation}
where $g_{j}=0$ for $j>p$ and $S_{0}=n$.

We now proceed case-by-case. The $p=2$ case illustrates the method,
so we will only include the intermediate steps in the $p=3,5,7$ cases.
\begin{caseenv}
\item $p=2$. Then, using (\ref{eq:Up-Def-W-Sum}), equation (\ref{eq:Up-Phi^t}) is
equivalent to\[
2^{8(n-1)}\big(2^{-1-12n}S_{n}\big)=2^{3}r,\text{ or}\]
\begin{equation}
S_{n}=2^{4n+12}r.\label{eq:St-Powers-of-2}\end{equation}
We now use (\ref{eq:Newtons-Formula}) to calculate $S_{1}$ and $S_{2}$:\[
S_{1}=g(1)\]
\[
S_{2}=g_{1}S_{1}-2g_{2}=g_{1}^{2}-2g_{2}.\]
From Lemma \ref{lem:Modular-Eq-Phi} we can compute the values of
the $g_{j}$. Using the $b_{j}$ values from the table in that lemma,
we have\[
g_{1}=2^{14}(b_{1}\phi_{2}+b_{2}\phi_{2}^{2})=2^{16}(3\phi_{2}+2^{8}\phi_{2}^{2})\]
\[
g_{2}=-2^{14}b_{2}\phi_{2}=-2^{24}\phi_{2}.\]
We can now see that \[
S_{1}=g_{1}=2^{16}(3\phi_{2}+2^{8}\phi_{2}^{2})\]
\[
S_{2}=2^{32}(3\phi_{2}+2^{8}\phi_{2}^{2})^{2}+2^{25}\phi_{2}=2^{20}(2^{5}\phi_{2}+2^{12}3^{2}\phi_{2}^{2}+2^{21}3\phi_{2}^{3}+2^{28}\phi_{2}^{4}).\]
Thus (\ref{eq:St-Powers-of-2}) is satisfied for $n=1,2$. We proceed
by induction. Assume (\ref{eq:St-Powers-of-2}) is satisfied for all
integers $<n$. We show that it is satisfied for $n$. For ease of
computation, we introduce the set \[
R^{*}=2^{8}R^{(2)}=\left\{ \sum_{i=1}^{m}d_{i}2^{8i}\phi_{2}^{i}\big|d_{i}\in\mathbb{Z},m\in\mathbb{Z}^{+}\right\} \]
which, the reader will notice, is a ring. From (\ref{eq:Newtons-Formula})
we obtain\[
S_{n}=g_{1}S_{n-1}-g_{2}S_{n-2}=2^{8}r_{1}^{*}\cdot2^{4n}r_{2}^{*}+2^{16}r_{3}^{*}\cdot2^{4(n-1)}r_{4}^{*}=2^{4n+8}r_{5}^{*}=2^{4n+16}r\]
where $r_{i}^{*}\in R^{*}$ and $r\in R^{(2)}$.
%Therefore, (\ref{eq:St-Powers-of-2})
%is satisfied for all integers $\geq1$, which completes this case.
\item $p=3$. We want to show \begin{equation}
S_{n}=3^{2n+7}r\label{eq:St-Powers-of-3}\end{equation}
 where $r\in R^{(3)}$. We compute the $g_{j}$ and $S_{n}$ as follows,
using the $b_{j}$ from the table:\[
\begin{array}{ccc}
g_{1}=3^{9}(3^{9}\phi_{3}^{3}+3^{5}4\phi_{3}^{2}+10\phi_{3}), &
g_{2}=3^{14}(-3^{4}\phi_{3}^{2}-4\phi_{3}), &
g_{3}=3^{18}\phi_{3},\end{array}\]
\[
\begin{array}{ccc}
S_{1}=g_{1}, & S_{2}=g_{1}^{2}-2g_{2}, &
S_{3}=g_{1}{}^{3}-3g_{1}g_{2}+3g_{3}.\end{array}\] From this, we
obtain\[
S_{1}=3^{9}(3^{9}\phi_{3}^{3}+3^{5}4\phi_{3}^{2}+10\phi_{3})\]
\[
S_{2}=3^{14}(8\phi_{3}+3^{5}34\phi_{3}^{2}+3^{9}80\phi_{3}^{3}+3^{13}68\phi_{3}^{4}+3^{18}8\phi_{3}^{5}+3^{25}\phi_{3}^{6})\]
\[
S_{3}=3^{19}(\phi_{3}+3^{5}40\phi_{3}^{2}+3^{8}1174\phi_{3}^{3}+3^{15}136\phi_{3}^{4}+3^{18}581\phi_{3}^{5}+3^{25}16\phi_{3}^{6}+3^{27}58\phi_{3}^{7}+3^{32}4\phi_{3}^{8}+3^{35}\phi_{3}^{9})\]
which proves (\ref{eq:St-Powers-of-3}) for $n=1,2,3$. For the inductive
step, let $R^{*}$ be the ring $3^{4}R^{(3)}$ so that\begin{align*}
S_{n} & =g_{1}S_{n-1}-g_{2}S_{n-2}+g_{3}S_{n-3}\\
 & =3^{5}r_{1}^{*}3^{2n+1}r_{2}^{*}+3^{10}r_{3}^{*}3^{2n-1}r_{4}^{*}+3^{14}r_{5}^{*}3^{2n-3}r_{6}^{*}\\
 & =3^{2n+6}r_{7}^{*}\\
 & =3^{2n+10}r\end{align*}
where $r_{i}^{*}\in R^{*}$ and $r\in R^{(3)}$.
%The $p=3$ case is
%complete.
\item $p=5$. We want\begin{equation}
S_{n}=5^{2n+2}r\label{eq:St-Powers-of-5}\end{equation}
where $r\in R^{(5)}$. Computing the $S_{n}$ we find\[
\begin{array}{ccc}
S_{1}=5^{5}r_{1} & S_{2}=5^{8}r_{2} & S_{3}=5^{10}r_{3}\end{array}\]
\[
\begin{array}{cc}
S_{4}=5^{13}r_{4} & S_{5}=5^{16}r_{5}\end{array}\]
for some $r_{1},\ldots,r_{5}\in R^{(5)}$. This proves (\ref{eq:St-Powers-of-5})
for $n=1,\ldots,5$. For the inductive step, let $R^{*}$ be the ring
$5R^{(5)}$ so that\begin{align*}
S_{n} & =g_{1}S_{n-1}-g_{2}S_{n-2}+g_{3}S_{n-3}-g_{4}S_{n-4}+g_{5}S_{n-5}\\
 & =5^{4}r_{1}^{*}5^{2n-1}r_{2}^{*}-\ldots+5^{14}r_{9}^{*}5^{2n-9}r_{10}^{*}\\
 & =5^{2n+3}r_{11}^{*}\\
 & =5^{2n+4}r\end{align*}
where $r_{i}^{*}\in R^{*}$ and $r\in R^{(5)}$ .
\item $p=7$. We want\begin{equation}
S_{n}=7^{n+2}r\label{eq:St-Powers-of-7}\end{equation}
where $r\in R^{(7)}$. Computing the $S_{n}$ we find\[
\begin{array}{cccc}
S_{1}=7^{4}r_{1} & S_{2}=7^{6}r_{2} & S_{3}=7^{7}r_{3} & S_{4}=7^{9}r_{4}\end{array}\]
\[
\begin{array}{ccc}
S_{5}=7^{11}r_{5} & S_{6}=7^{13}r_{6} & S_{7}=7^{15}r_{7}\end{array}\]
for some $r_{1},\ldots,r_{7}\in R^{(7)}$. This proves (\ref{eq:St-Powers-of-7})
for $n=1,\ldots,7$. For the inductive step, let $R^{*}$ be the ring
$7R^{(7)}$ so that\begin{align*}
S_{n} & =\sum\limits _{i=1}^{7}(-1)^{i+1}g_{i}S_{n-i}\\
 & =7^{3}r_{1}^{*}7^{n}r_{2}^{*}-\ldots+7^{13}r_{13}^{*}7^{n-6}r_{14}^{*}\\
 & =7^{n+3}r_{15}^{*}\\
 & =7^{n+4}r\end{align*}
where $r_{i}^{*}\in R^{*}$ and $r\in R^{(7)}$ .
\end{caseenv}
\end{proof}

\section{Proof of the Theorem}

To remind the reader of the main result of the paper, we include it
here.
\begin{thm*}
Let $p\in\{2,3,5,7\}$, and let $f_{0,m}^{(p)}(\tau)=q^{-m}+\sum
a_{0}^{(p)}(m,n)q^{n}$ be an element of the basis described above,
with $m=p^{\alpha}m'$ and $(m',p)=1$. Then, for $\beta > \alpha$,
\[
\begin{array}{rll}
a_{0}^{(2)}(2^{\alpha}m',2^{\beta}n)\equiv0 & \pmod{2^{3(\beta-\alpha)+8}} & \text{if }p=2\\
a_{0}^{(3)}(3^{\alpha}m',3^{\beta}n)\equiv0 & \pmod{3^{2(\beta-\alpha)+3}} & \text{if }p=3\\
a_{0}^{(5)}(5^{\alpha}m',5^{\beta}n)\equiv0 & \pmod{5^{(\beta-\alpha)+1}} & \text{if }p=5\\
a_{0}^{(7)}(7^{\alpha}m',7^{\beta}n)\equiv0 & \pmod{7^{(\beta-\alpha)}} & \text{if }p=7.\end{array}\]
\end{thm*}
The proof is in three cases. The first illustrates the method for
the simplest basis elements, namely those with $(m,p)=1$. The second
demonstrates the ``shifting'' property at its first occurence,
$f_{0,p}^{(p)}$. The third is the general case; it builds
inductively upon the methods of the first two cases.

\subsection{Case 1: $(m,p)=1$}
\begin{proof}
This proof is almost identical to Lehner's proof of Theorem 3 in
\cite{Lehner:2}, however it applies not only to functions which
have poles of order bounded by $p$, but to all basis elements with
$(m,p)=1$. For ease of notation, let $f(\tau)=f_{0,m}^{(p)}(\tau)$.

We will demonstrate the method with $m=1$, then generalize it to
all $m$ relatively prime to $p$. First, we will write $U_{p}f(\tau)$
as a polynomial in $\phi(\tau)$ with integral coefficients, all of
which are divisible by the desired power of $p$. Since $U_{p}$ isolates
the coefficients whose index is divisible by $p$, we will have proven
the theorem for $\beta=1$. We will then apply $U_{p}$ repeatedly
to the polynomial, showing that the result is always another polynomial
in $\phi$ with integral coefficients, all of which are divisible
by the desired higher power of $p$.

Consider the level $p$ modular function $g(\tau)=pU_{p}f(\tau)+p^{\lambda/2}\phi(\tau)$.
Notice that $g(\tau)$ is holomorphic at $\infty$ since both $U_{p}f(\tau)$
and $\phi(\tau)$ are holomorphic there. The $q$-expansion at $0$
for $g(\tau)$ is given by\[
g(-1/p\tau)=p(U_{p}f)(-1/p\tau)+p^{\lambda/2}\phi(-1/p\tau)\]
which, by Lemmas \ref{lem:Psi-At-0} and \ref{lem:Main-Up-Formula} becomes\[
g(-1/p\tau)=p(U_{p}f)(p\tau)+f(-1/p^{2}\tau)-f(\tau)+\psi(\tau).\]
When we notice that $f(\tau)=\psi(\tau)$ in this $m=1$ case, we
obtain\begin{align*}
g(-1/p\tau) & =p(U_{p}f)(p\tau)+\psi(-1/p^{2}\tau)-\psi(\tau)+\psi(\tau)\\
 & =p(U_{p}f)(p\tau)+p^{\lambda/2}\phi(p\tau),\end{align*}
which is holomorphic at $\infty$. Hence, $g(\tau)$ is a holomorphic
modular function on $\Gamma_{0}(p)$, so it must be constant. Therefore,
\begin{equation}
U_{p}f(\tau)=c_{0}-p^{\lambda/2-1}\phi(\tau)\label{eq:Up-f1=00003DPhi}\end{equation}
for some constant $c_{0}$. The proof is complete for $\beta=1$.

Note: the prime $13$, having genus zero, would work in this construction;
however, in that case $\lambda=\frac{24}{13-1}=2$, so $13^{\lambda/2-1}=1$,
and we gain no new information.

We now iterate the above process to prove the theorem for $\beta>1$.
Notice that \[
U_{p}(U_{p}f(\tau))=c^{(p)}-p^{\lambda/2-1}U_{p}\phi(\tau).\] We
know from Lemma \ref{lem:Modular-Eq-Phi} that $U_{p}\phi$ is a
polynomial in $\phi$; in fact, by inspection of the $b_{j}^{(p)}$
values we see that we may write\begin{align*}
U_{2}\phi^{(2)}(\tau) & =2^{3}\big(d_{1}^{(2)}\phi^{(2)}(\tau)+\sum_{n=2}^{2}d_{n}^{(2)}2^{8(n-1)}\phi^{(2)}(\tau)^{n}\big)\\
U_{3}\phi^{(3)}(\tau) & =3^{2}\big(d_{1}^{(3)}\phi^{(3)}(\tau)+\sum_{n=2}^{3}d_{n}^{(3)}3^{4(n-1)}\phi^{(3)}(\tau)^{n}\big)\\
U_{5}\phi^{(5)}(\tau) & =5\big(d_{1}^{(5)}\phi^{(5)}(\tau)+\sum_{n=2}^{5}d_{n}^{(5)}5^{n}\phi^{(5)}(\tau)^{n}\big)\\
U_{7}\phi^{(7)}(\tau) &
=7\big(d_{1}^{(7)}\phi^{(7)}(\tau)+\sum_{n=2}^{7}d_{n}^{(7)}7^{n}\phi^{(7)}(\tau)^{n}\big)\end{align*}
for some integers $d_{n}^{(p)}$. This shows that the second $U_{p}$
iteration is divisible by the correct power of $p$. Further, it
gives us a polynomial of a suitable form to iterate the process
using Lemma \ref{lem:Phi-Polynomials}. In each of the polynomials
above, notice that $U_{p}\phi(\tau)=p^{\delta}r$ for some $r\in
R^{(p)}$. Using Lemma \ref{lem:Phi-Polynomials}, we find that\[
U_{p}(U_{p}\phi)(\tau)=p^{2\delta}r'\] for some $r'\in R^{(p)}$, and
further\[ U_{p}^{\beta}\phi(\tau)=p^{\beta\delta}r^{(\beta)}\] for
some $r_{\beta}\in R^{(p)}$. This completes the proof for $m=1$.

Now, if $(m,p)=1$, then $U_{p}f(\tau)$ is holomorphic at $\infty$,
just as it was with $m=1$. Moving to the cusp at $0$ we find
that $(U_{p}f)(-1/p\tau)$ can be written as a polynomial in $\psi(\tau)$
which appears as a polynomial in $\phi(\tau)$ when we return to $\infty$.
Similar to (\ref{eq:Up-f1=00003DPhi}), we obain the equality\[
U_{p}f(\tau)=c_{0}+\sum_{i=1}^{M}p^{\lambda i/2-1}c_{i}\phi(\tau)^{i}\]
for some $c_{i}\in\mathbb{Z}$ and $M\in\mathbb{Z}^{+}$. The only
difference between this equation and (\ref{eq:Up-f1=00003DPhi}) is
that in this more general case, we find that $U_{p}f$ is a higher
degree polynomial in $\phi$. This formula can easily be iterated
as before to obtain the desired result.
\end{proof}

\subsection{Case 2: $m=p$}
\begin{proof}
Again, for ease of notation, denote $f_{0,p}^{(p)}(\tau)$ by
$f(\tau)$. For the $m=p$ case, we will proceed as before; however,
we will find that $U_{p}f(\tau)$ has poles at both $\infty$ and $0$,
and that $U_{p}f(\tau)$ does not possess any interesting
divisibility properties, but $U_{p}^{2}f(\tau)$ does. This property
will manifest itself as the {}``shifting'' previously mentioned.

Notice first that $U_{p}f(\tau)=q^{-1}+O(1)$ has a simple pole at
$\infty$. Therefore, we shall deal primarily with the function $U_{p}f(\tau)-\psi(\tau)$,
which is holomorphic at $\infty$. We can use Lemmas \ref{lem:Psi-At-0}
and \ref{lem:Main-Up-Formula} to view this function at $0$:

\begin{align*}
p(U_{p}f)(-1/p\tau)-p\psi(-1/p\tau) & = p(U_{p}f)(p\tau)+f(-1/p^{2}\tau)-f(\tau)-p^{\lambda/2+1}\phi(\tau) \\
 & =pq^{-p}+O(1)+O(1)-q^{-p}+O(1)+O(q)\\
 & =c_{0}+\sum_{i=1}^{p}c_{i}\psi(\tau)^{i}
\end{align*}
for some integers $c_{i}$. Replacing $\tau$ by $-1/p\tau$, we obtain
\begin{equation}
(U_{p}f)(\tau)=\frac{c_{0}}{p}+\psi(\tau)+\sum_{i=1}^{p}c_{i}p^{\lambda
i/2-1}\phi(\tau)^{i}.\label{eq:Case-2-Upf}
\end{equation}

The $\psi(\tau)$ term in the equation makes any attempt at
$p$-divisibility fail; for example, computation shows that the
$7^{th}$ coefficient of $\psi^{(2)}(\tau)$ is odd. However,
$\psi(\tau)$ satisfies Lehner's divisibility properties, so $U_{p}f$
inherits its $p$-divisibility from $\psi(\tau)$. So the function\[
U_{p}^{2}f(\tau)=c_{0}+U_{p}\psi(\tau)+\sum_{i=1}^{p}c_{i}p^{\lambda
i/2-1}U_{p}\phi(\tau)^{i}\] has the same $p$-divisibility as
$f_{0,1}^{(p)}$; hence, the shift.
\end{proof}

\subsection{Case 3: $m=p^{\alpha}m'$}
\begin{proof}
We prove this case using induction on $\alpha$. Case 1 showed that the theorem is true for
all $m'$ relatively prime to $p$, so the $\alpha=0$ base case is complete. Assume Theorem
\ref{thm:Andersen-Main}
holds for all $m$ of the form $m=p^{\ell}m'$ with $\ell<\alpha$.
We will show it holds for $m=p^{\alpha}m'$.
To simplify notation,
let $f_{\alpha}(\tau)=f_{0,p^{\alpha}m'}^{(p)}(\tau)$.

Since $f_{\alpha}(\tau)=q^{-p^{\alpha}m'}+O(1)$, we find that
$U_{p}f_{\alpha}(\tau)=q^{-p^{\alpha-1}m'}+O(1)$ has a pole of order
$p^{\alpha-1}m'$ at $\infty$. So we focus our attention on
$U_{p}f_{\alpha}(\tau)-f_{\alpha-1}(\tau)$, which is holomorphic at
$\infty$. Using (\ref{eq:Main-Up-Formula}) we examine this function
at $0$:\begin{align*}
p(U_{p}f_{\alpha})\left(\frac{-1}{p\tau}\right) -
pf_{\alpha-1}\left(\frac{-1}{p\tau}\right) &
=p(U_{p}f_{\alpha})(p\tau)+f_{\alpha}\left(\frac{-1}{p^{2}\tau}\right)
-
f_{\alpha}(\tau)-pf_{\alpha-1}\left(\frac{-1}{p\tau}\right)\\
 & =pq^{-p^{\alpha}m'}+O(1)+O(1)-q^{-p^{\alpha}m'}-O(1)-O(1)\\
 & =(p-1)q^{-p^{\alpha}m'}+O(1).\end{align*}

As before, we write this function as a polynomial in $\psi(\tau)$
with integral coefficients $c_{i}$:\[
p(U_{p}f_{\alpha})(-1/p\tau)-pf_{\alpha-1}(-1/p\tau)=c_{0}+\sum_{i=1}^{p^{\alpha}m'}c_{i}\psi(\tau)^{i},\]
which, after switching back to the $q$-expansion at $\infty$, becomes\begin{equation}
U_{p}f_{\alpha}(\tau)=\frac{c_{0}}{p}+f_{\alpha-1}(\tau)+\frac{1}{p}\sum_{i=1}^{p^{\alpha}m'}c_{i}p^{\lambda i/2}\phi(\tau)^{i}.\label{eq:Case-3-Upf}\end{equation}

Notice that (\ref{eq:Case-3-Upf}) looks very similar to
(\ref{eq:Case-2-Upf}), where $\psi(\tau)$ is replaced by
$f_{\alpha-1}(\tau)$, so $U_{p}f_{\alpha}(\tau)$ inherits whatever
divisibility properties $f_{\alpha-1}(\tau)$ has. Our inductive
hypothesis states that $f_{\alpha-1}(\tau)$ exhibits Lehner's
divisibility properties only after $U_{p}$ is applied $\alpha-1$
times. Therefore, applying $U_{p}$ to (\ref{eq:Case-3-Upf})
$\alpha-1$ times, we obtain\[
U_{p}^{\alpha}f_{\alpha}(\tau)=\frac{c_{0}}{p}+U_{p}^{\alpha-1}f_{\alpha-1}(\tau)
+\frac{1}{p}\sum_{i=1}^{p^{\alpha}m'}c_{i}p^{\lambda
i/2}U_{p}^{\alpha-1}\phi(\tau)^{i}\] showing that
$U_{p}^{\alpha}f_{\alpha}(\tau)$ exhibits Lehner's divisibility
properties.
\end{proof}
\bibliographystyle{amsplain}
%\bibliography{bibliography}

\begin{thebibliography}{1}

\bibitem{Ahlgren:theta}
Scott Ahlgren.
\newblock The theta-operator and the divisors of modular forms on genus zero
  subgroups.
\newblock {\em Math. Res. Lett.}, 10(5-6):787--798, 2003.

\bibitem{Apostol:modular}
Tom~M. Apostol.
\newblock {\em Modular functions and {D}irichlet series in number theory},
  volume~41 of {\em Graduate Texts in Mathematics}.
\newblock Springer-Verlag, New York, second edition, 1990.

\bibitem{DoudJenkinsLopez}
D.~Doud, P.~Jenkins, and J.~Lopez.
\newblock Two-divisibility of the coefficients of certain weakly holomorphic
  modular forms.
\newblock arXiv:1105.0684v1 [math.NT].

\bibitem{Doud:padic}
Darrin Doud and Paul Jenkins.
\newblock {$p$}-adic properties of coefficients of weakly holomorphic modular
  forms.
\newblock {\em Int. Math. Res. Not. IMRN}, (16):3184--3206, 2010.

\bibitem{Duke:zeros}
W.~Duke and Paul Jenkins.
\newblock On the zeros and coefficients of certain weakly holomorphic modular
  forms.
\newblock {\em Pure Appl. Math. Q.}, 4(4, Special Issue: In honor of
  Jean-Pierre Serre. Part 1):1327--1340, 2008.

\bibitem{Griffin}
M.~Griffin.
\newblock Divisibility properties of coefficients of weight 0 weakly
  holomorphic modular forms.
\newblock to appear in International Journal of Number Theory, DOI
  10.1142/S1793042111004599.

\bibitem{Lehner:1}
Joseph Lehner.
\newblock Divisibility properties of the fourier coefficients of the modular
  invariant $j(\tau)$.
\newblock {\em Amer. J. Math.}, 71(1):136--148, Jan 1949.

\bibitem{Lehner:2}
Joseph Lehner.
\newblock Further congruence properties of the fourier coefficients of the
  modular invariant $j(\tau)$.
\newblock {\em Amer. J. Math.}, 71(2):373--386, Apr 1949.

\end{thebibliography}

\end{document}